\documentclass[12pt,twoside]{article}
\usepackage{amssymb,latexsym}
\newtheorem{theorem}{Theorem}
\newtheorem{propo}[theorem]{Proposition}
\newtheorem{lem}[theorem]{Lemma}
\newtheorem{cor}[theorem]{Corollary}

\begin{document}\vspace*{0pt}

\def\sgn{{\rm sgn}}

\let\a=\alpha
\let\d=\delta
\let\D=\Delta
\let\e=\varepsilon
\let\f=\varphi
\let\g=\gamma
\let\G=\Gamma
\let\i=\iota
\let\l=\lambda
\let\L=\Lambda
\let\m=\mu
\let\n=\nu
\let\N=\nabla
\let\nd=\nabla
\let\p=\pi
\let\r=\rho
\let\s=\sigma
\let\S=\Sigma
\let\u=\theta
\let\w=\omega
\let\W=\Omega
\let\x=\xi
\let\X=\Xi
\let\y=\psi

\def\bbR{{\mathbb{R}}}
\def\bbZ{{\mathbb{Z}}}
\def\bbC{{\mathbb{C}}}
\def\cA{{\cal A}}
\def\cP{{\cal P}}

\def\jth{j^{\underline{{\rm th}}}}
\newcommand{\usfrac}[2]{{#1}/{#2}}
\def\endo{{\rm End}}

\newcommand{\dfrac}[2]{{\displaystyle{\frac{#1}{#2}}}}
\newcommand{\tfrac}[2]{{\textstyle{\frac{#1}{#2}}}}
\newcommand{\nn}[1]{(\ref{#1})}
\newcommand{\nnn}[1]{{\rm(}\ref{#1}{\rm)}}

\begin{center}{\Large\bf Spontaneous generation of eigenvalues}

{\large\sf Thomas Branson and Bent {\O}rsted}
\end{center}

\begin{center}\begin{minipage}{14cm}
{\bf Abstract}: We show that the action of conformal vector fields on 
functions on the 
sphere determines the spectrum of the Laplacian 
(or the conformal Laplacian), without
further input of information.  
The spectra of intertwining operators (both differential and non-local)
with principal part a power of the Laplacian follows as a corollary.
An application of the method is the sharp form of Gross' entropy
inequality on the sphere.  The same method gives the spectrum of the 
Dirac operator
on the sphere, as well as of a continuous family of 
nonlocal intertwinors, and an infinite family of odd-order differential
intertwinors.
\end{minipage}\end{center}

\subsection*{Introduction}
Among compact Riemannian manifolds, the sphere is remarkable in
several ways. For example, it admits a large group of conformal transformations,
and the corresponding conformal vector fields are known to
contain important information about the geometry of the sphere.
What is perhaps not so well-known is the fact that the action of
these vector
fields also determines a large part of global analysis on the
sphere, notably the exact eigenvalues of most natural differential 
and pseudodifferential operators.  Already the
example of the Laplace operator on functions is perhaps a little surprising.

In \cite{spg}, a method was developed for finding the spectrum of
intertwining operators for certain representations of semisimple
groups.  In this calculation scheme, one first inputs the spectral
data of a differential 
{\em spectrum generating operator}.  The idea is that this
spectral data is readily accessible, obtained just from quadratic
Casimir data for finite-dimensional representations of compact groups.
Using this information, one may generate the (much less accessible)
spectral data on the intertwining operators, which occur in series
parameterized by a complex number $r$.  Typically, a subseries of {\em
differential} operators occurs at values of $r$ which are positive
integral, in a sense appropriate to the particular series.

A special situation that may occur in such a calculation is when
the spectrum generating operator is closely related to one of
the differential 
intertwinors.  The simplest interpretations of ``closely related''
here are (1) when the operators are identical, or differ by a 
constant additive shift; (2) when the spectrum generating operator is
a polynomial in one of the differential intertwinors.  
In some of these 
cases, there may be some many relations among the spectral data that
one can generate the spectra of the operators in question
with no input at all; that is, the spectra are {\em spontaneously
generated}.  

In this paper, we execute this process of spontaneous generation 
in some examples, and derive some related results.  Among the central
relations are a quadratic equation \nn{quadr} relating ``adjacent''
eigenvalues of the conformal Laplacian on scalar fields, and a cubic equation
\nn{cubic} which serves a similar purpose for the Dirac operator on spinor
fields.
Once the spectra of these fundamental operators are in place, one may
go on to find parameterized families of operators of various orders
which are {\em intertwining} for conformal transformations and vector
fields in the same sense as are the conformal Laplacian and Dirac operator.

An interesting feature of this construction is that it proves completeness
of the eigenvalue list it obtains -- any ``wrong'' eigenvalue will generate
lower and lower ones, until some basic estimate is violated.  And 
even though no information on spherical harmonics is input, 
the fact that eigenfunctions are spherical harmonics is a consequence
of the construction.
(See, for example, Corollary \ref{protosh}.)

Our method, combined with Beckner's sharp Hardy-Littlewood-
Sobolev inequalities on the sphere, is perfectly suited for deriving
the sharp form of Gross' entropy inequality on the sphere,
estimating the integral of $f^2\log(f)$ for a positive smooth function.
This argument culminates in Theorem \ref{entro} below.

Of course, some of the spectral resolutions we get here have been
known for some time; for example, the spectra of the Laplacian
(on functions), and of the Dirac operator on the sphere.  The spectra
of the other operators is less well known, but nevertheless
can be obtained by
specializing formulas in the literature (for example \cite{spg}).
However, the philosophy is avoid all heavy machinery, and derive
the spectra just from a couple of elementary operator commutation
relations.  Among the possible virtues of this sort of ``primitivist''
derivation is the prospect of deriving such spectra at an early
point of a course in differential geometry, quantum mechanics, or
relativity.

\subsection*{The Laplacian and conformal Laplacian on the sphere}

Let $n\ge 2$, and let
$\r_0$ be the {\em azimuthal angle} on the sphere $S^n$ with homogeneous
coordinate
functions $x_0,x_1,\cdots,x_n$.  That is, 
$$
x_0=\cos\r_0,
$$
so that $\r(p)$ is the angle between an indeterminate point $p$ and the point
$(1,0,\cdots,0)$.  More invariantly, we can define an azimuthal angle
$\r(p,q)$
with any desired point $q$ in place of $(1,0,\cdots,0)$.  In 
particular, the azimuthal angle from the point where $x_i=1$ will
be denoted $\r_i$ (for each $i$).

The proper conformal
vector fields on $S^n$ are generated by the $T_i:=
(\sin\r_i)(\partial/\partial\r_i)$ for $j=0,\cdots,n$.
The meaning of $\partial/\partial\r_i$ is: complete
$\r_i$ to a coordinate system by 
taking a coordinate system on
the latitude $\{\r_i={\rm const}\}$ (a copy of the sphere $S^{n-1}$), 
and compute with these other coordinates
held constant.  
The $T_i$ are conformal vector fields:
$$
L(T_i)g=2x_ig,
$$
where $g$ is the round metric on $S^n$ and 
$L$ denotes the Lie derivative.  

In the more general setting of pseudo-Riemannian manifolds $(M,g)$, let
$d$ be the exterior derivative on functions, and let $\d$ be its formal
adjoint.
The
operator $\D=\d d$ is the {\em Laplacian}, and the {\em conformal Laplacian}
is the operator
$$
D:=\D+\frac{n-2}{4(n-1)}K,
$$
where $K$ is the scalar curvature of $g$.  This satisfies the {\em conformal 
covariance} relation
\begin{equation}\label{ccy}
\W^{(n+2)/2}D_{\W^2g}f=D_g(\W^{(n-2)/2}f)  
\end{equation}
for $\W>0$ and $f$ smooth functions on $M$.
If $(M,g)$ admits a conformal vector field $T$ with $L_Tg=2\w g$, 
then \nn{ccy} immediately leads,
via the local flow, to
\begin{equation}\label{infcc}
D\left(T+\frac{n-2}2m(\w)\right)=\left(T+\frac{n+2}2m(\w)\right)D,
\end{equation}
where $m(\w)$ is multiplication by $\w$.
On the other hand,
\begin{equation}\label{protospecgen}
[D,m(\f)]=[\D,m(\f)]=m(\D f)-2\i(d\f)d,
\end{equation}
where $\i$ is interior multiplication.

Specializing to the round sphere $S^n$, let
$$
U_j:=T_j+\frac{n}2x_j.
$$
Then \nn{infcc} reads
\begin{equation}\label{ccready}
D(U_i-m(x_i))=(U_i+m(x_i))D.
\end{equation}
In particular, applying this relation to the function 1 and using
the fact that the scalar curvature of round $S^n$ is $n(n-1)$, so that
\begin{equation}\label{roundconflap}
D=\D+\frac{n(n-2)}4\ \ \mbox{ on }S^n,
\end{equation}
we get
$$
\left(\D+\frac{n(n-2)}4\right)\frac{n-2}2x_i
=\frac{n+2}2\cdot x_i\cdot\frac{n(n-2)}4,
$$
that is,
\begin{equation}\label{foothold}
\D x_i=nx_i.
\end{equation}

Specializing \nn{protospecgen}, 
we have
\begin{equation}\label{specgen}
\begin{array}{rl}
[D,m(x_i)]&=m(\D x_i)-2\i(dx_i)d \\
&= nm(x_i)+2T_i \\
&=2U_i,
\end{array}
\end{equation}
where we have used \nn{foothold} together with the fact that
$-dx_i$ corresponds to $T_i$ under the metric
identification.

\newcommand{\pres}[3]{{}_{#1}|\,{#2}|_{#3}}
Let $E(\l)$ denote the $\l$-eigenspace of $D$, and suppose for
some $\l$ we have $0\ne\f\in E(\l)$.
Given another number $\m$ and an operator $T$, denote the compression of $T$
acting $E(\l)\to E(\m)$ by $\pres{\m}{T}{\l}$.

Given $\m$, if we compress \nn{specgen} and \nn{ccready} to operators
from $E(\l)$ to $E(\m)$, we have
\begin{equation}\label{compress}
\begin{array}{rl}
\pres{\m}{U_i}{\l}&=
\dfrac{\m-\l}2\pres{\m}{m(x_i)}{\l}, \\
\m\cdot\pres{\m}{(U_i-m(x_i))}{\l}&=\l\cdot\pres{\m}{(U_i+m(x_i))}{\l},
\end{array}
\end{equation}
from which
$$
\m\left(\frac{\m-\l}2-1\right)\pres{\m}{m(x_i)}{\l}=
\l\left(\frac{\m-\l}2+1\right)\pres{\m}{m(x_i)}{\l}.
$$
For fixed $\l$, this is a quadratic equation in $\m$:
\begin{equation}\label{quadr}
\m^2-2\l\m-2\m+\l^2-2\l=0.
\end{equation}
The solutions 
\begin{equation}\label{process}
\l^\pm
:=\l+1\pm\sqrt{4\l+1};
\end{equation}
are candidates for new eigenvalues of $D$.
Since $\D$ is a nonnegative operator, \nn{roundconflap} gives
\begin{equation}\label{esti}
\l\ge\frac{n(n-2)}4.
\end{equation}
Thus 
$$
\l^-<\l<\l^+
$$
unless $n=2$ and $\l=0$, in which case $\l^-=\l$.

Moreover, since \nn{quadr} is symmetric in $\m$ and $\l$,
we have that $\m$ is a root of the equation obtained starting with $\l$
iff $\l$ is a root of the equation obtained starting with $\m$.
That is,
$$
(\l^+)^-=\l=(\l^-)^+.
$$

The obvious main questions are:
\begin{itemize}
\item If $E(\l)\ne 0$, are the spaces $E(\l^\pm)$ necessarily nonzero?
\item Can we generate all $\m$ for which $E(\m)\ne 0$ by starting
with a single $E(\l)\ne 0$ and iterating the process \nn{process}?
\end{itemize}
The answer to each question is yes, with the exception that 
$E(\l^-)$ vanishes for what should be the bottom eigenvalue $\l$,
namely $n(n-2)/4$.  In fact, both questions are answered by the same
calculation: any ``wrong'' eigenvalue keeps generating lower and
lower ``wrong'' eigenvalues, until we get one that violates the
estimate \nn{esti}.

A weakness of the calculation is that it cannot immediately tell
the dimension of the eigenspaces $E(\l)$.  Of course these dimensions
are easily obtainable with an injection of a small amount of Lie
theory, but that would violate the spirit of the present calculation.
Note however, that in Corollary 6 we prove (from our point of view)
the correspondence with spherical harmonics.

Stipulating the statements just above (which will be proved below),
we may start either
with the fact that
$$
1\in E(0,\D)=E(n(n-2)/4,D)
$$
or 
$$
x_i\in E(n,\D)=E((n+2)n/4,D)\ \ \ \mbox{(from \nn{foothold})}
$$
and iterate \nn{process}
to compute that the eigenvalues of $D$ are the
\begin{equation}\label{theeigs}
\l_j:=
\left(\frac{n-2}2+j\right)\left(\frac{n}2+j\right)\ \mbox{ for }\ j=0,1,2,\cdots,
\end{equation}
or equivalently, that the eigenvalues of $\D$ are the $j(n-1+j)$.

To begin to answer the bullet point questions above, let us look
more closely at the eigenfunctions (as opposed to just eigenvalues)
that our process is generating.  If $\f\in E(\l,D)$, then
\begin{equation}\label{whatReigfcns}
D(x_i\f)=\l x_i\f+2U_i\f
\end{equation}
by \nn{specgen}.
By \nn{ccready} and \nn{whatReigfcns},
\begin{equation}\label{whatReigtwo}
D(U_i\f)=D(x_i\f)+\l(U_i+x_i)\f
=(\l+2)U_i\f+2\l x_i\f.
\end{equation}
If we would like $(U_i+cx_i)\f$
to be in $E(\m,D)$, then
$$
\m(U_i+cx_i)\f=(\l+2+2c)U_i\f+(2+c)\l x_i\f.
$$
If the $U_i$ and $x_i$ terms on the two sides are to agree,
we need 
$$
c(\l+2+2c)=\l(2+c).
$$
This is a quadratic equation on $c$, with the roots
$$
c^\pm=
(\l-\l^\mp)/2.
$$
The eigenvalue $\m$ corresponding to $c^+$ (resp.\ $c^-$) is 
easily computed to be $\l^+$ (resp.\ $\l^-$).  Thus we have:

\begin{propo}\label{PandM} If $\f\in E(\l,D)$, then
$$
\begin{array}{rl}
P_i\f:=&\left(U_i+\frac12(\l-\l^-)
x_i\right)
\f\in E(\l^+,D), \\
M_i\f:=&\left(U_i+\frac12(\l-\l^+)
x_i
\right)\f\in E(\l^-,D).
\end{array}
$$
\end{propo}

Of course there is {\em a priori} the danger that $P_i\f$ vanishes
for $i=0,\cdots,n$,
or that each $M_i\f$ does, even though $\f\ne 0$.  We now proceed to rule
out this danger, except in the case where 
it is expected (applying $M_i$ to the eigenspace where $1$ lives).
That is, we shall compute
$$
\sum_iM_iP_i\f\ \mbox{ and }\ \sum_iP_iM_i\f.
$$
If the first of these is nonzero, it must be that $P_i\f\ne 0$ for some
$i$; similarly, if the second is nonzero, some $M_i\f$ does not vanish
identically.
We begin with
$$
\sum_iM_iP_i\f=\sum_i\left(U_i+\frac{-1-\sqrt{1+4\l^+}}2x_i\right)
\left(U_i+\frac{-1+\sqrt{1+4\l}}2x_i\right)\f.
$$
Since
\begin{equation}\label{sqrts}
\sqrt{1+4\l^\pm}=\pm 2+\sqrt{1+4\l},
\end{equation}
the above is
$$
\begin{array}{l}
\displaystyle\sum_i\left(U_i+\frac{-3-\sqrt{1+4\l}}2x_i\right)
\left(U_i+\frac{-1+\sqrt{1+4\l}}2x_i\right)\f
= \\
\displaystyle\sum_i\left((U_i-m(x_i))+c^-m(x_i)\right)
\left((U_i-m(x_i))-c^-m(x_i)\right)\f.
\end{array}
$$
This simplifies to
\begin{equation}\label{midway}
\begin{array}{l}
\Bigg\{\sum_i(U_i-m(x_i))^2-(c^-)^2m\left(\sum_ix_i^2\right) \\
-c^-\left(\sum_i\left\{(U_i-m(x_i))m(x_i)-m(x_i)(U_i-m(x_i))\right\}\right)
\Bigg\}\f.
\end{array}
\end{equation}
Now 
\begin{equation}\label{sumone}
\sum_ix_i^2=1,
\end{equation}
and commuting this relation with $D$,
\begin{equation}\label{anticomm}
\begin{array}{rl}
0&=[D,m(1)] \\
&=\sum_i\{m(x_i)[D,m(x_i)]+[D,m(x_i)]m(x_i)] \\
&=2\sum_i\{m(x_i)U_i+U_im(x_i)\}
\end{array}
\end{equation}
by \nn{specgen}.  This shows in turn that
\begin{equation}\label{orth}
\sum_i(U_i+am(x_i))^2=a^2+\sum_iU_i^2
\end{equation}
for any number $a$.  On the other hand,
\begin{equation}\label{comm}
\begin{array}{rl}
\sum_i[U_i,m(x_i)]&=\sum_iT_i\cos\r_i=-\sum_i\sin^2\r_i \\
&=-\sum_i(1-\cos^2\r_i)=1-(n+1)=-n.
\end{array}
\end{equation}

It is clear that we shall also need a simplification
of $\sum_iU_i^2$.  This is provided by:

\begin{lem}\label{TandD} $\sum_iT_i^2=-\D$.
\end{lem}

{\em Proof:} Consider one term of the sum on the left, and suppress the subscript
$i$ for now, so that $T_i=T$, $x_i=x$.  In abstract index notation, one term
from the 
left side of the identity is
\begin{equation}\label{ab}
x^a\nd_a(x^b\nd_b)=(\nd^ax)\{(\nd^bx)\nd_a+\nd_a\nd^bx\}\nd_b.
\end{equation}
Since $\nd^ax=(dx)^a$, the first term is
$(dx\otimes dx)^{ab}\nd_a\nd_b$.  The second term involves
$$
\nd_a\nd_bx=({\rm Hess}\,x)_{ba}.
$$
But the conformal Killing equation $L_Tg=2xg$ reads, in abstract index notation,
$$
\nd_aT_b+\nd_bT_a=2xg_{ab}.
$$
Since $T_a=-(dx)_a$, this says that $2$Hess$\,x=-2xg$, so that the quantity \nn{ab} becomes
$$
(dx\otimes dx)^{ab}\nd_a\nd_b-(\nd^ax)x\d_a{}^b\nd_b
=(dx\otimes dx)^{ab}\nd_a\nd_b-(\nd^ax)x\nd_a.
$$
Re-inserting the subscript $i$ and summing over it,
the above quantity becomes
$$
g^{ab}\nd_a\nd_b-\frac12(\nd^a1)\nd_a=g^{ab}\nd_a\nd_b=-\D,
$$
since $\sum_idx_i\otimes dx_i$ is the pullback of the ambient flat 
$\bbR^{n+1}$ metric, i.e.\ the round metric, and
$$
\sum_ix_i\nd^ax_i=\frac12\nd^a\underbrace{\left(\sum_ix_i^2\right)}_1.\qquad\square
$$

By the lemma and \nn{orth},
$$
-\left(D-\dfrac{n(n-2)}4\right)=
\sum_iU_i^2+\dfrac{n^2}4,
$$
so that
$$
\sum_iU_i^2=-D-n/2.
$$
Using all the identities just derived to evaluate \nn{midway}, we get
$$
\sum_iM_iP_i\f=\left(-\l-\frac{n}2+1-(c^-)^2+nc^-\right)\f.
$$
It is convenient to write this in terms of
$$
\n:=\sqrt{1+4\l},
$$
so that
$$
\l=\dfrac{\n^2-1}4,\qquad c^-=-\dfrac{1+\n}2,
$$
and
\begin{equation}\label{MPgets}
\sum_iM_iP_i\f=-\frac12(\n+n-1)(\n+2)\f.
\end{equation}
This shows that some $P_i\f$ is nonzero unless $\n=1-n$ or $\n=-2$;
in particular it is nonzero for all positive values of $\n$.
But by \nn{esti}, $\n\ge n-1$.

The corresponding calculation with $M_i$ and $P_i$ in the other
order begins with
$$
\sum_iP_iM_i\f=\sum_i\left(U_i+\frac{-1+\sqrt{1+4\l^-}}2x_i\right)
\left(U_i+\frac{-1-\sqrt{1+4\l}}2x_i\right)\f,
$$
and (by virtue of \nn{sqrts}) produces a version of \nn{midway}
with $c^+$ in place of $c^-$, yielding
$$
\sum_iP_iM_i\f=\left(-\l-\frac{n}2+1-(c^+)^2+nc^+\right)\f.
$$
Since $c^+=(\n-1)/2$, we get
\begin{equation}\label{PMgets}
\sum_iP_iM_i\f=-\frac12(\n-n+1)(\n-2)\f.
\end{equation}
Thus some $M_i\f$ is nonzero unless $\n=n-1$ or $\n=2$.
In the first case, $\l=n(n-2)/4$; this was expected,
since \nn{esti} and $1\in E(n(n-2)/4,D)$ show that this is the 
bottom eigenvalue of $D$.
Indeed, a look back at the formula for $M_i\f$ in Proposition \ref{PandM} shows
that $M_i1=0$.  In the second case, $\n=2$, we have $\l=3/4$.
By \nn{esti}, this implies that $n=3$ and equality holds in \nn{esti};
this is a special case of the situation just discussed.  We have proved:

\begin{propo}\label{getthere} For $\l\ge n(n-2)/4$, if $0\ne\f\in E(\l,D)$,
then $P_i\f$ is a nonzero element of $E(\l^+,D)$ for some $i$.
If $\l>n(n-2)/4$, then 
$M_k\f$ is a nonzero element of $E(\l^-,D)$ for some $k$.
\end{propo}

The mechanism by which we have generated the eigenvalues
\nn{theeigs} also rules out any other numbers occurring as eigenvalues.
Indeed, suppose $\l$ is an eigenvalue of $D$ not on the list, so
that (by \nn{esti}) there is some natural number $j$ with
$\l_j<\l<\l_{j+1}$.
Since the map $\l\to\l^-$ is strictly monotonic, we have
$\l_{j-1}<\l^-<\l_j$, $\l_{j-2}<\l^{--}<\l_{j-1}$ and so on,
until we reach eigenvalues $\m\in(\l_0,\l_1)$, $\m^-<\l_0$,
contradicting \nn{esti}.  The key ingredient, of course, is
the assurance from Proposition \ref{getthere} that each lower eigenspace in
the induction is truly nonzero.  We have:

\begin{propo}\label{complete} The $\l_j$ of {\rm\nn{theeigs}} give the complete list
of eigenvalues of $D$ on $S^n$.  As a result, the $j(n-1+j)$
for $j=0,1,2,\cdots$ give the complete list of eigenvalues of $\D$
on $S^n$.
\end{propo}

We can also harvest the following corollary.  Let $E_j:=E(\l_j,D)$.

\begin{cor}\label{xgetsthere} 
The span of the $m(x_i)E_j$ and $Dm(x_i)E_j$ is
$E_{j-1}\oplus E_{j+1}$
(where $E_{-1}=0$).
\end{cor}

{\em Proof}: The inclusion $\subset$ is immediate from the first line of 
\nn{compress} together with Proposition \ref{PandM}.  For the inclusion
$\supset$, (\ref{MPgets},\ref{PMgets}) show that $E_{j-1}\oplus E_{j+1}$
is contained in the sum of the $m(x_i)E_j$ and the 
$$
U_iE_j\subset Dm(x_i)E_j+m(x_i)DE_j=Dm(x_i)E_j+m(x_i)E_j.
\qquad\square
$$

Note that we have not injected any information on spherical harmonics into
our procedure for generating the eigenvalues.  We may, however, get the
interpretation of the eigenspaces as spaces of spherical harmonics as
a consequence of what we have done:

\begin{cor}\label{protosh}
$E_j$ is exactly the set of restrictions from $\bbR^{n+1}$ 
to $S^n$ of $j$-homogeneous harmonic polynomials in 
the $x_i$.
\end{cor}

{\em Proof}: First, the elements of $E_j$ are restrictions
of $j$-homogeneous polynomials, since this is true of $E_0$, so follows from
the previous corollary by induction on $j$.
Second, if we compute the Laplacian of $\bbR^{n+1}$ in spherical coordinates,
we get
$$
\D_{\bbR^{n+1}}=-\dfrac{\partial^2}{\partial r^2}
-\dfrac{n-1}{r}\dfrac{\partial}{\partial r}+\dfrac1{r^2}\D_{S^n},
$$
so that if $\f\in E_j$, and we extend to $\bbR^{n+1}$ by extending the 
$j$-homogeneous polynomial formula, 
$$
\D_{\bbR^{n+1}}\f=[-j^2-(n-1)j+j(n-1+j)]\f=0.
$$
This calculation
also shows that each $j$-homogeneous harmonic polynomial in $\bbR^{n+1}$ gives rise
to an element of $E_j.\qquad\square$

As a bonus result, we can give operators $A_{2r}$ which are functions of $D$
(or of $\D$), which satisfy generalizations of the conformal covariance
relation \nn{ccready}, namely
\begin{equation}\label{ccr}
D(U_i-r\cdot m(x_i))=(U_i+r\cdot m(x_i))D,
\end{equation}
for each $r\in\bbC$, each of which takes an eigenvalue on $E_j:=E(\l_j,D)$.
That is, we can find functions $f_r$ on the spectrum of $D$ for which
$$
A_{2r}|_{E_j}=f_r(\l_j){\rm Id}_{E_j}.
$$
In the usual notation of functional calculus, we write $A_{2r}=f_r(D)$.

{\bf Remark}: Though we have shown 
that we have all eigenvalues for 
the conformal Laplacian, our derivation does not contain a proof that
the span of the corresponding eigenfunctions is dense in $L^2(S^n)$.
Of course, we have this by general elliptic theory.  This sort of 
completeness is implicitly used later, when we describe other
covariant operators, or {\em intertwinors}, as functions of a basic
one (for example $D$, or the operator $A_1$ of \nn{aT}, or 
the Dirac operator $P$ as it is used in Theorem \ref{spinortwinors}.

First note that $\l_j^\pm=\l_{j\pm 1}$.  Compressing
\nn{ccr}, and letting $\m_j$ be the putative eigenvalue for $A_{2r}$ on 
$E_j$, we have
$$
\m_{j\pm 1}\left(\dfrac{\l_{j\pm 1}-\l_j}2-r\right)\pres{E_{j\pm 1}}{m(x_i)}{E_j}=
\m_j\left(\dfrac{\l_{j\pm 1}-\l_j}2+r\right)\pres{E_{j\pm 1}}{m(x_i)}{E_j}.
$$
Since the $\l_j$ are known, we may generate the various $\m_j$ inductively
by demanding
\begin{equation}\label{recur}
\m_{j\pm 1}\left(\l_{j\pm 1}-\l_j-2r\right)=
\m_j\left(\l_{j\pm 1}-\l_j+2r\right).
\end{equation}
Since {\em a priori} this gives two relations between adjacent $\m_j$,
we must check for consistency.
Since $\l_{j+1}-\l_j=n+2j$, we have this, provided we handle
occurrences of vanishing $\l_{j\pm 1}-\l_j\pm'2r$ correctly.
We emerge with a choice of $\m_j^{(2r)}$ that is unique up to a 
constant (independent
of $j$) factor:

\begin{propo}\label{bonus} For fixed $r\notin\{-n/2,-n/2-1,\cdots\}$,
\begin{equation}\label{specfcn}
Z(r,j):=\dfrac{\G(n/2+j+r)}{\G(n/2+j-r)}
\end{equation}
is, up to a constant nonzero factor, 
the unique function of $j\in\{0,1,2,\cdots\}$ that satisfies \nnn{recur} and 
does not vanish identically.  For $r=-n/2-j_0$ with $j_0$ a nonnegative integer, 
the residue of the above expression {\rm(}viewed as a meromorphic function of $r${\rm)}
at $-n/2-j$ is, up to a constant nonzero factor, 
the unique function of $j\in\{0,1,2,\cdots\}$ that satisfies \nnn{recur} and 
does not vanish identically.
\end{propo}

If $r\in\frac12\bbZ^+$, then 
$$
Z(r,j)=(n/2+j+r-1)\cdots(n/2+j-r).
$$
If $-r\in\frac12\bbZ^+$, then
$$
Z(r,j)=\left(\displaystyle
\prod_{1\le p\le-2r}^{\bullet}
(n/2+j-r-p)\right)^{-1},
$$
where $\displaystyle\prod^\bullet$ is the product over nonzero factors.

The operator $A_1$ has eigenvalue $Z(\frac12,j)=\usfrac{(n-1)}2+j$ on $E_j$;
thus
\begin{equation}\label{aT}
A_1:=\sqrt{\D+\left(\frac{n-1}2\right)^2}.
\end{equation}
This implies via \nn{specfcn} that
$$
A_{2r}=\dfrac{\G(A_1+\frac12+r)}{\G(A_1+\frac12-r)},\ \ r\notin\{-n/2,-n/2-1,\cdots\}.
$$
Similarly, $A_{2r}$ (for $r$ outside the exceptional set given above) may be
written as a function of $A_{2q}$ for any nonzero $q$ outside the exceptional set.
In particular, for $r\in\bbZ^+$ we get the sequence of differential operators
$$
\displaystyle\prod_{p=1}^r\left\{\D+\left(\frac{n}2+p-1\right)
\left(\frac{n}2-p\right)\right\},
$$
also written down in \cite{tbjfa}, Remark 2.23.

\subsection*{The entropy inequality}

Equation \nn{recur}, together with the sharp Hardy-Littlewood-Sobolev
inequalities of Beckner \cite{bec} on the sphere, give an argument for
an optimal form of Gross' entropy inequality on the sphere.

Let $r$ be the parameter of \nn{recur}, and denote the $r$-derivative
at $r=0$ by a prime.  Differentiating \nn{recur} and normalizing
so that the intertwinor $A_0$ is the identity, we have
$$
\m'_{j\pm 1}(\l_{j\pm 1}-\l_j)-2
=\m'_j(\l_{j\pm 1}-\l_j)+2.
$$
Thus
\begin{equation}\label{info}
(\m_{j\pm 1}-\m_j)'=\dfrac{4}{\l_{j\pm 1}-\l_j}\,.
\end{equation}
The $j+1$ and $j-1$ relation lists are consistent, so the information
in \nn{info} is equivalent to
\begin{equation}\label{infoH}
(\m_{j+1}-\m_j)'=\dfrac{4}{\l_{j+1}-\l_j}\,.
\end{equation}
Since 
$$
\l_{j+1}-\l_j=n+2j,
$$
this gives (with $m:=n/2$)
$$
\m'_j=\m'_0+\dfrac2{m}+\dfrac2{m+1}+\cdots+\dfrac2{m+j-1}\,.
$$
This is in fact a formula
for the eigenvalues of $A'_{2r}\,$.

It is convenient to pick a normalization of the series $A_{2r}$ with
the property that $\m'_0=0$.  This is obtained by 
multiplying the spectral function \nn{specfcn} by $\G(m-r)/\G(m+r)$
-- this factor is independent of $j$, so for fixed $r$, gives a constant
multiple of the $A_{2r}$ described by \nn{specfcn}.  With this
normalization,
$$
\m'_j=\dfrac2{m}+\dfrac2{m+1}+\cdots+\dfrac2{m+j-1}\,.
$$
Let us denote the intertwinors normalized in this way by $B_{2r}$.
These operators also appear in Beckner's sharp
Hardy-Littlewood-Sobolev
inequalities
for $r\in[0,n/2)$: if $F$ is a positive smooth function on the sphere,
$$
\int F^{(n-2r)/2}B_{2r}F^{(n-2r)/2}\ge\left(\int F^n\right)^{(n-2r)/n},
$$
where all integrals are with respect to normalized measure.
One has equality for every $F$ when $r=0$ (each side is $\int F^n$).
But for $r>0$, equality holds exactly when $F$ is a constant multiple
of a conformal (diffeomorphism) factor.

Now write the inequality as
$$
0\le
-\left(\int F^n\right)^{(n-2r)/n}
+\int F^{(n-2r)/2}B_{2r}F^{(n-2r)/2},
$$
and take $(d/dr)|_{r=0}$ of each side.
(Of course this differentiation of the inequality depends 
on the fact that equality holds for every $F$ at $r=0$.)
This gives
$$
0\le\dfrac2{n}\left(\int F^n\right)\log\int F^n
-2\int F^n\log F+\int F^{n/2}B'_{2r}F^{n/2},
$$
or in a slightly better form,
$$
2\int F^n\log F\le\dfrac2{n}\left(\int F^n\right)\log\int F^n
+\int F^{n/2}B'_{2r}F^{n/2}.
$$
What we know immediately about the case of equality is that it includes 
{\em at least} constant multiples of 
the conformal factors $F$.  But since the quantities in
play are analytic in $r$ near $r=0$, we get {\em precisely}
these functions.

With $f:=F^{n/2}$, we can rewrite as
$$
\dfrac4{n}\int f^2\log f\le\dfrac2{n}\left(\int f^2\right)\log\int f^2
+\int fA'_{2r}f,
$$
for
$$
{\rm eig}(A'_{2r},E_j)=
\dfrac2{m}+\dfrac2{m+1}+\cdots+\dfrac2{m+j-1}.
$$
It is easily verified that these eigenvalues are $\le$ (but 
very close to) those of $2\log(2A_1/(n-1))$, where $A_1$ is as in
\nn{aT}.  (The harmonic sum is a certain 
Riemann sum for the integral defining
the log.)
Thus we may write
\begin{equation}\label{giveaway}
\dfrac2{n}\int f^2\log f\le\dfrac1{n}\left(\int f^2\right)\log\int f^2
+\int f\left(\log\dfrac{2A_1}{n-1}\right)f,
\end{equation}
giving away a little sharpness.
To summarize:

\begin{theorem}\label{entro} For smooth positive $f$ on the sphere $S^n$, in
normalized
measure,
\begin{equation}\label{lHLS}
\dfrac4{n}\int f^2\log f\le\dfrac2{n}\left(\int f^2\right)\log\int f^2
+\int fA'_{2r}f,
\end{equation}
where $H$ takes the eigenvalue 
$$
\dfrac2{m}+\dfrac2{m+1}+\cdots+\dfrac2{m+j-1}
$$
on $\jth$ order spherical harmonics.  
In particular the weaker statement \nnn{giveaway} holds.
Equality holds in 
\nnn{lHLS} if and only if $f^{2/n}$ is a positive constant multiple
of a conformal {\rm(}diffeomorphism{\rm)} factor.
\end{theorem}

\subsection*{Spinor operators}

Let $P$ be the Dirac operator on the spinor bundle $\Sigma$.
By the Lichnerowicz
formula,
$$
P^2=\N^*\N+\frac14K,
$$
so that 
$$
P^2=\N^*\N+\frac14n(n-1)\ \mbox{ on }S^n.
$$
As a result, if $\l$ is an eigenvalue of $P$ on the sphere, then
$\l$ is real with
\begin{equation}\label{estiS}
\l^2\ge n(n-1)/4.
\end{equation}

The analogue of \nn{specgen} is
\begin{equation}\label{diracspg}
[P^2,m(\w)]=[\nd^*\nd,m(\w)]=2\nd_T+nm(\w).
\end{equation}
The conformal covariance relation satisfied by $P$ is
\begin{equation}\label{interdirac}
P\left(L_T+\frac{n-1}2\w\right)=\left(L_T+\frac{n+1}2\w\right)P,
\end{equation}
on general pseudo-Riemannian spin manifolds, where $L_Tg=2\w g$.
If $X$ is an arbitrary smooth vector field, \cite{kosm} shows that
the Lie and covariant
derivatives on spinors are related by
$$
L_X-\nd_X=-\frac18(dX)_{ab}\g^a\g^b,
$$
where $dX$ is the exterior derivative of the 1-form corresponding to $X$ under
the metric, and $\g$ is the fundamental tensor-spinor (a section of 
$TM\otimes\endo(\S)$).
This assumes that the internal conformal weight 0 has been assigned to
the spinor bundle; assigning internal weight $\pm\frac12$ as in
\cite{PR} results in an extra term involving div$\,X$.
Our proper conformal vector fields $T_i$ on $S^n$ have the
$-dx_i$ as their metric correspondents, so we may
specialize 
\nn{diracspg} to
\begin{equation}\label{diracspgtwo}
[P^2,m(x_i)]=2L(T_i)+nm(x_i)=:2U_i.
\end{equation}
The conformal covariance relation
\nn{interdirac} specializes to 
\begin{equation}\label{ccP}
P(U_i-\tfrac12m(x_i))=(U_i+\tfrac12m(x_i))P.
\end{equation}

Compressing the operators in 
\nn{diracspgtwo} to act between two eigenspaces for $P$, $E(\l,P)\to E(\m,P)$, we have
\begin{equation}\label{presU}
\pres{\m}{U_i}{\l}=\dfrac{\m^2-\l^2}2\pres{\m}{m(x_i)}{\l},
\end{equation}
after which \nn{ccP} implies that
$$
\m(\m^2-\l^2-1)\pres{\m}{m(x_i)}{\l}=\l(\m^2-\l^2+1)\pres{\m}{m(x_i)}{\l}.
$$
This is implied by
$$
\m(\m^2-\l^2-1)=\l(\m^2-\l^2+1),
$$
which may be rewritten
\begin{equation}\label{cubic}
(\m+\l)(\m-\l+1)(\m-\l-1)=0.
\end{equation}
(Note that all these equations are symmetric in $\m$ and $\l$.)

The cubic equation \nn{cubic} suggests the possible nonvanishing of three ``adjacent''
eigenspaces $E(-\l,P)$, $E(\l+1,P)$, and $E(\l-1,P)$, given $E(\l,P)\ne 0$.
Let $0\ne\y\in E(\l,P)$, and write the Dirac operator as $\g^a\nd_a$
(using abstract index notation, in which repetition of an index, once up and
once down, denotes a contraction).  To avoid excessive super- and subscripting,
denote $m(x_i)$ and $U_i$ by $x$ and $U$ for now.  In analogy with 
(\ref{whatReigfcns},\ref{whatReigtwo}), we have:
$$
\begin{array}{rl}
P(x\y)&=\g^a\nd_a(x\y)=\l x\y+\g^a(\nd_a x)\y, \\
P(U\y)&=P(\tfrac12x\y)+(U+\tfrac12x)\l\y
=\l x\y+\l U\y+\tfrac12\g^a(\nd_a x)\y.
\end{array}
$$
In contrast to (\ref{whatReigfcns},\ref{whatReigtwo}), however, this system doesn't
close.  We need in addition:
$$
P(\g^a(\nd_a x)\y)=\g^b\nd_b(\g^a(\nd_a x)\y)
=\g^b\g^a((\nd_a x)\nd_b\y+(\nd_b\nd_ax)\y),
$$
where we have used the relation $\nd\g=0$ between the spin connection and the
fundamental tensor-spinor.  By the Clifford relation
$$
\g^a\g^b+\g^b\g^a=-2g^{ab},
$$
we may rewrite the above as
$$
P(\g^a(\nd_a x)\y)=-\g^a(\nd_ax)\l\y-2(\nd^ax)\nd_a\y-(\nd^a\nd_ax)\y
=2U\y-\l\g^a(\nd_ax)\y,
$$
where we have used \nn{foothold} to simplify $-\nd^a\nd_ax=\D x$.

Abbreviating 
\begin{equation}\label{defy}
[P,m(x)]=\g^a(\nd_ax)=:y,
\end{equation}
we have
\begin{equation}\label{transit}
\begin{array}{rl}
A\y&:=(U+\l x+\tfrac12y)\y\in E(\l+1,P), \\
S\y&:=(U-\l x-\tfrac12y)\y\in E(\l-1,P), \\
N\y&:=(U-\tfrac12x-\l y)\y\in E(-\l,P).
\end{array}
\end{equation}
Recall that we have suppressed the subscript $i\in\{0,\cdots,n\}$, so that
we really have operators $A_i$, $S_i$, and $N_i$.  To find whether we really
get something nonzero by the above processes, we compute
\begin{equation}\label{saasnn}
\begin{array}{l}
\sum_iS_iA_i\y=\sum_i(U_i-(\l+1)x_i-\tfrac12y_i)(U_i+\l x_i+\tfrac12y_i)\y, \\
\sum_iA_iS_i\y=\sum_i(U_i+(\l-1)x_i+\tfrac12y_i)(U_i-\l x_i-\tfrac12y_i)\y, \\
\sum_iN_iN_i\y=\sum_i(U_i-\tfrac12x_i+\l y_i)(U_i-\tfrac12x_i-\l y_i)\y.
\end{array}
\end{equation}

In simplifying this, note that
$$
yy=\g^a\g^b(\nd_ax)\nd_bx=-(\nd_ax)\nd^ax=-|dx|^2.
$$
Thus
$$
\sum_iy_i^2=-\sum_i|dx_i|^2=-\sum_i\sin^2\r_i=-\sum_i(1-\cos^2\r_i)=1-(n+1)=-n.
$$
We also have
$$
yx=xy=x\g^a(\nd_ax)=\frac12\g^a\nd_a(xx),
$$
so that 
$$
\sum_iy_ix_i=\sum_ix_iy_i=\frac12\g^a\nd_a\left(\sum_ix_i^2\right)=\frac12\g^a\nd_a1=0.
$$
\nn{anticomm}, \nn{orth}, and \nn{comm} are still valid (with the new meaning of $U_i$,
and $P^2$ in place of $D$ in the intermediate steps of \nn{anticomm}).

In addition,
$$
\begin{array}{l}
Uy=U[P,x]=\frac12[P^2,x][P,x]=\frac12(P[P,x]+[P,x]P)[P,x]=Py^2+yPy, \\
yU=[P,x]U=\frac12[P,x][P^2,x]=\frac12[P,x](P[P,x]+[P,x]P)=yPy+y^2P,
\end{array}
$$
so that
$$
[U,y]=[P,y^2],
$$
and
\begin{equation}\label{Uy}
\sum_i[U_i,y_i]=\left[P,\sum_iy_i^2\right]=[P,-n]=0.
\end{equation}

The argument of Lemma \ref{TandD} goes through formally as written
(with $\nd$ now involving the spin connection), and gives
$$
\sum_i\nd_{T_i}^2=-\nd^*\nd=-P^2+\dfrac{n(n-1)}4.
$$
With the discussion preceding \nn{diracspgtwo}, this gives
$$
\sum_i\left(U_i-\tfrac{n}2x_i\right)^2=-P^2+\dfrac{n(n-1)}4.
$$
By the analogues of (\ref{anticomm},\ref{orth},\ref{comm}), 
$$
\begin{array}{rl}
\sum_i(U_i+ax_i)^2&=a^2+\sum_iU_i^2 \\
&=a^2-P^2+\dfrac{(n(n-1)}4-\dfrac{n^2}4 \\
&=a^2-P^2-\dfrac{n}4.
\end{array}
$$
In particular,
$$
\sum_iU_i^2=-P^2-\dfrac{n}4.
$$

Using all these identities, for arbitrary $a$ and $b$,
$$
\sum_i(U_i-(a+1)x_i-by_i)(U_i+ax_i+by_i)=-P^2-\dfrac{3n}4-a^2-a(n+1)+b^2n.
$$
In the first line of \nn{saasnn}, we take $a=\l$ and $b=\frac12$, and apply
to a $\l$-eigenspinor $\y$ to get
$$
\sum_iS_iA_i\y=\left(-2\l^2-\frac{n}2-(n+1)\l\right)\y=-2\left(\l+\frac{n}2\right)
\left(\l+\frac12\right)\y.
$$
In the second line, take $a=-\l$ and $b=-\frac12$ to get
$$
\sum_iA_iS_i\y=\left(-2\l^2-\frac{n}2+(n+1)\l\right)\y=-2\left(\l-\frac{n}2\right)
\left(\l-\frac12\right)\y.
$$
In the third line, take $a=-\frac12$ and $b=-\l$ to get
$$
\sum_iN_iN_i\y=(n-1)\left(\l+\frac{1}2\right)\left(\l-\frac12\right)\y.
$$
This establishes:

\begin{propo} If $P$ has a nonzero eigenspace on $S^n$, then 
the nonzero eigenspaces of $P$ on round $S^n$ are exactly
the 
$$
F_j:=E(n/2+j,P),\qquad G_j:=E(-(n/2+j),P)
$$
for $j=0,1,2,\cdots$.  
\end{propo}

{\em Proof:} The last three displayed identities and \nn{estiS} show that if $E(\l,P)\ne 0$,
then 
$$
\begin{array}{l}
E(-\l,P)\ne 0, \\
E(\l-1,P)\ne 0\qquad{\rm unless}\ \l=n/2, \\
E(\l+1,P)\ne 0\qquad{\rm unless}\ \l=-n/2.
\end{array}
$$
Thus if there is an eigenvalue outside the set of $\pm(n/2+j)$, we may generate
an eigenvalue contradicting \nn{estiS}.  Given any nonzero eigenspace $E(\l,P)$
with $\l$ of the form $\pm(n/2+j)$ however, 
we may generate nonzero eigenspaces 
corresponding to {\em all} such $\l.\qquad\square$

The beginning of the statement of the last proposition is a bit awkward;
we need to assume there {\em is} some eigenspace in order to get a foothold
analogous to that provided by the function 1 in the scalar case.  We can
get this foothold by taking a nonzero parallel spinor in $\bbR^n$ and 
stereographically injecting it (with the proper conformal weight)
to $S^n$.  Alternatively, we could do just enough elementary elliptic theory
to conclude that the minimizer for the Rayleigh quotient based on $P^2$
provides us with an eigensection.

The analogue of Corollary \ref{xgetsthere} is 

\begin{cor}\label{xygets}
The span of the $m(x_i)E(\l,P)$, $Pm(x_i)E(\l,P)$, and $P^2m(x_i)E(\l,P)$
is $E(\l+1,P)\oplus E(\l-1,P)\oplus E(-\l,P)$.
\end{cor}

{\em Proof:} The inclusion $\subset$ follows from \nn{transit}, 
\nn{presU}, and the fact (from \nn{defy}, in the notation of that
display) that $\pres{\m}{y}{\l}=(\m-\l)\pres{\m}{m(x)}{\l}$.
For the inclusion $\supset$, note first that \nn{saasnn}
puts $E(\l+1,P)\oplus E(\l-1,P)\oplus E(-\l,P)$
in the span of the $A_iE(\l,P)$, $N_iE(\l,P)$, and $S_iE(\l,P)$.
By \nn{saasnn}, this is in the span of the 
$x_iE(\l,P)$, $U_iE(\l,P)$, and $y_iE(\l,P)$.
By \nn{diracspgtwo} and \nn{defy}, this is in the span of the 
$m(x_i)E(\l,P)$, $Pm(x_i)E(\l,P)$, and $P^2m(x_i)E(\l,P).\qquad
\square$

In analogy with Proposition \ref{bonus}, we may seek intertwinors
$\cA_{2k+1}$ for complex-valued $k$, satisfying the intertwining relation
\begin{equation}\label{kint}
\cA_{2k+1}(U-(k+\tfrac12)x)
=(U+(k+\tfrac12)x)\cA_{2k+1}
\end{equation}
(which extends \nn{ccP}).  In fact, such intertwinors exist and
commute
with $P$.  The subfamily of these 
for which $k$ is a nonnegative 
integer yields a family of odd-order (in fact, order
$2k+1$) differential operators which are polynomial in $P$:

\begin{theorem}\label{oddpoly} For $k=0,1,2,\cdots$, and
$$
\cP_{2k+1}=(P-k)(P-k+1)\cdots P\cdots(P+k-1)(P+k)
=P(P^2-1)(P^2-4)\cdots(P^2-k^2),
$$
we have
$$
\cP_{2k+1}(U-(k+\tfrac12)x)
=(U+(k+\tfrac12)x)\cP_{2k+1},
$$
where $U$ is any $U_i$, and $x$ is any $x_i$.
\end{theorem}

{\em Proof}: Fix $k$, and consider 
$$
\a(\l):=(\l-k)(\l-k+1)\cdots\l\cdots(\l+k-1)(\l+k),
$$
the eigenvalue 
taken by $\cP_{2k+1}$ on 
$E(\l,P)$.
In view of \nn{presU}, it is enough to show that
\begin{equation}\label{ets}
\begin{array}{l}
(\m-k)(\m-k+1)\cdots\l\cdots(\m+k-1)(\m+k)\left(\dfrac{\m^2-\l^2}2
-(k+\tfrac12)\right)\stackrel{?}{=} \\
\qquad(\l-k)(\l-k+1)\cdots\l\cdots(\l+k-1)(\l+k)
\left(\dfrac{\m^2-\l^2}2
+(k+\tfrac12)\right)
\end{array}
\end{equation}
for $\m=\l\pm 1$ and for $\m=-\l$.
This is easily verified in each of the three cases.$\qquad\square$

For arbitrary order $k+1/2$, we might expect a nonlocal intertwining
operator.  This is provided by:

\begin{theorem}\label{spinortwinors} If $k+n/2\notin\bbZ$, 
the operators
\begin{equation}\label{pp}
\cA_{2k+1}:=
\sgn(P)^{n+1}\,\dfrac{\Gamma(P+k+1)}{\G(P-k)}
\end{equation}
satisfy the intertwining relations \nn{kint}, where for
each function $f$, the operator $f(P)$ takes the eigenvalue $f(\l)$
on $E(\l,P)$.
\end{theorem}

{\em Proof}: Fix $k$ and 
let $\a[\l]$ be the eigenvalue on $E(\l,P)$ of a putative 
intertwinor $\cA_{2k+1}$.  It is immediate from \nn{presU} that
$\a[-\l]=-\a[\l]$, and that 
\begin{equation}\label{shift}
\a[\l+1](\l-k)=(\l+1+k)\a[\l].
\end{equation}
Moreover, these equations for all $P$-eigenvalues $\l$ are sufficient
for the intertwining property (the condition relating $\a[\l-1]$
to $\a[\l]$ is the same as \nn{shift}, with $\l$ shifted to $\l-1$.)
Implementing the recursion this gives, we get the intertwinor
\begin{equation}\label{spinorint}
\sgn(P)\,\dfrac{\Gamma(|P|+k+1)}{\G(|P|-k)},
\end{equation}
and one may check directly that this is intertwining.
Using the identity
$$
\G(z)\G(-z)z\sin(\pi z)=-\pi,
$$
we get
$$
\dfrac{\G(-P+k+1)}{\G(-P-k)}=\dfrac{\G(P+k+1)}{\G(P-k)}
\cdot\dfrac{\sin\pi(P+k)}{\sin\pi(P-k-1)}.
$$
The trigonometric factor in this is
$$
\dfrac{\sin(\pi P)\cos(\pi k)+\cos(\pi P)\sin(\pi k)}
{\sin(\pi(P-1))\cos(\pi k)-\cos(\pi(P-1))\sin(\pi k)}\,.
$$
If $n$ is even, $k\notin\bbZ$ 
and $P$ takes integral eigenvalues, so this becomes
$1$.  If $n$ is odd, $P$ takes properly half-integral eigenvalues
and $k$ is not a proper half-integer, so the above becomes $-1$.
This shows that 
$$
\dfrac{\Gamma(|P|+k+1)}{\G(|P|-k)}=\sgn(P)^n
\dfrac{\Gamma(P+k+1)}{\G(P-k)},
$$
so that the intertwinor \nn{spinorint} is the $\cA_{2k+1}$ of 
\nn{pp}.$\qquad\square$

{\bf Remark}: When $n$ is even, we may set $k=1/2$ in the last theorem
to obtain the intertwinor
\begin{equation}\label{weird}
\sgn(P)
\dfrac{\Gamma(P+\tfrac32)}{\G(P-\tfrac12)}=\sgn(P)(P+\tfrac12)(P-\tfrac12)
=P|P|-\tfrac14\usfrac{P}{|P|}.
\end{equation}
This is also a formula for an intertwinor when $n$ is odd, though this is not
covered by Theorem \ref{spinortwinors}.  To see this, we just feed the formula
for the operator \nn{weird} into a calculation similar to \nn{ets}.  Similar
considerations hold for the other values of $k$ excluded by the hypotheses of
Theorem \ref{spinortwinors}.  Of course it is known that even and odd
dimensions act very differently in many ways with respect to the spinor
bundle.  For example, in even dimensions the bundle $\S$ itself admits a
chirality decomposition into two subbundles.  In contrast, $\S$ is irreducible
in odd dimensions, but each eigenspace of $P^2$ admits a chiral decomposition.
In particular, the relation between $E(\l,P)$ and $E(-\l,P)$ has very
different interpretations in the two dimension parities.  The technique we use
here seems remarkably attuned to deriving Theorem \ref{spinortwinors}
using only things that are common to the two dimension parities, the 
power on $\sgn(P)$
being the only hint of the differences.

{\bf Remark}: Our technique can also be used to obtain the eigenvalue specta
of $D$, $P$ and other intertwinors on the non-Riemannian hyperboloids
\cite{rs}; i.e.\ on other real forms of the sphere, again manifolds with large
Lie algebras of conformal vector fields.  This provides some explanation, in
terms of conformal geometry, of the remarkable agreement between most of the
eigenvalues of the Laplace operator on the various real forms of the sphere.

Thomas Branson: Department of Mathematics, University of Iowa, Iowa City
IA 52242 USA\hfill
{\tt thomas-branson@uiowa.edu}

Bent {\O}rsted: Institut for Matematiske Fag, Aarhus Universitet,
8000 Aarhus C, Denmark\hfill
{\tt orsted@imf.au.dk}

\end{document}